\documentclass[12pt,reqno]{amsart}
\usepackage{amsmath,amsfonts,amssymb,amsthm,hyperref,enumerate,multicol}
\usepackage{graphicx}
\usepackage{tikz}
\usetikzlibrary{shapes.geometric, arrows}
\newtheorem{theorem}{Theorem}

\newtheorem{example}{Example}

\numberwithin{equation}{section}
\begin{document}

\title[]
{On wavelet-based sampling Kantorovich operators and their study in multi-resolution analysis}

\maketitle

\begin{center}
{\bf Digvijay Singh$^1$ Rahul Shukla$^2$, Karunesh Kumar Singh$^3$ } \footnote{Corresponding author: Karunesh Kumar Singh} \\
\vskip0.15in
$^{1,3}$ Department of Applied Sciences and Humanities, Institute of Engineering and Technology, Lucknow, 226021, Uttar Pradesh, India \\
$^{2}$ Department of Mathematics, Deshbandhu College, University of Delhi, India \\
\vskip0.15in

\vskip0.15in

Email: dsiet.singh@gmail.com$^1$, rshukla@db.du.ac.in$^2$,   kksiitr.singh@gmail.com$^3$,
\end{center}

\begin{abstract}
In this work, wavelet-based filtering operators are constructed by introducing a basic function $D(t_1, t_2, t_3)$ using a general wavelet transform. The cardinal orthogonal scaling functions (COSF) provide an idea to derive the standard sampling theorem in multiresolution spaces which motivates us to study wavelet approximation analysis. With the help of modulus of continuity, we establish a fundamental theorem of approximation. Moreover, we unfold some other aspects in the form of an upper bound of the estimation taken between the operators and functions with various conditions. In that order, a rate of convergence corresponding to the wavelet-based filtering operators is derived, by which we are able to draw some important interferences regarding the error near the sharp edges and smooth areas of the function.
Eventually, some examples are demonstrated and empirically proven to justify the fact about the rate of convergence. Besides that, some derivation of inequalities with justifications through examples and important remarks emphasizes the depth and significance of our work. \\

\end{abstract}


\section{Introduction}
In the late 1980s, the wavelet theory was coined by Mallat \& Daubechies \cite{I1}. Roughly speaking, wavelet theory is a mathematical model that analyses small waves to represent signals. From the application point of view, wavelets theory has been one of the most growing areas in mathematics, especially in wave propagation, data compression, image processing, pattern recognition, computer graphics, the detection of aircrafts, and image compression etc. 
\par In wavelets theory, one of the main tasks is to find the coefficients of a polynomial using the inner product space, and sometimes it becomes very difficult because of the rapid oscillation as the frequency goes to infinity and oscillates very slowly at a very low frequency. For the first time, Morlet and his significant group of researchers \cite{I1} face this sort of difficulty and have handled this hurdle efficiently by reconstructing the suitable function. Afterwards, in 1998, Valery A. Zheludev \cite{Z1} developed a generalized form of wavelets called `mother wavelets' which is used to produce a collection of wavelets which is defined as follows
$$ \Upsilon_{a_1,a_2}(t)=\frac{1}{\sqrt{|a_1|}}\Upsilon\left(\frac{t-a_2}{a_1}\right)\,\, a_1,a_2\in \mathbb{R}\,\text{with}\, a_1\,\neq 0.$$
Moreover, in order to study the wavelets theory deeply by introducing in form of dilation and scaling factors, we have a term called `wavelets transform' developed by a group of researchers including Jean Morlet, Yves Meyer, Stéphane Mallat, and Ingrid Daubechies \cite{A2,I11,A3,A1}. Particularly, if a function $f\in L^2(\mathbb{R})$ then the general continuous wavelet transform \cite{I2,I3} is defined as follows
\begin{eqnarray}
    \label{W1}
    (\mathcal{W}_{\Upsilon} f)(a_1,a_2)=\int_{\mathbb{R}} f(x) \overline{\Upsilon_{a_1,a_2}(x)} dx,
\end{eqnarray}
where $a_1$ and $a_2$ are called scale and time parameters, respectively. The notation $\overline{}$ denotes the complex conjugate.
If we fix the scale time plane \((a_1, a_2)\) at certain levels of \(a_1 = 2^{-j}\), \(a_2 = \varrho 2^{-J}\), where \(j < J\) and \(\varrho \in \mathbb{Z}\), for some \(J > 0\), the wavelet transform denoted by \((\mathcal{W}_{\Upsilon} f)(a_1, a_2)\) is given as follows;

\begin{eqnarray*}
    (\mathcal{W}_{\Upsilon} f)(a_1, a_2)=WCT\{f(t); a_1=2^{-j}\,, a_2=\varrho2^{-j}\}=b_{j,\varrho}\,\,,(say)\, j,\,\varrho\,\in \mathbb{Z}.
\end{eqnarray*}
Mallat and Shensa \cite{I11,I10} developed a useful algorithm by which, the representation of signals at a specific time scale became feasible. The calculation of coefficients \( b_{J,\varrho} \) contains some tedious integrals. Therefore, in order to solve it, there has always been a leverage using these algorithms.\par
From the above discussion, one can think what would be the representation of discrete to continuous wavelets? A rational approach tends us to think that it depends on the nature of the signals. Therefore, the classical Shannon sampling theorem is used to reconstruct the function completely to a band-limited in a compact interval in such a way that the rate of uniform samples exceeds the Nyquist rate. The results based on the classical fundamental theorem have numerous applications in the area of signal processing and communication theory, can be seen in \cite{I5,I4}. Moreover, the idea has been extended by Walter Aldroubi, and Unser \cite{I7,I6} to multi-resolution spaces comes from the area called `multi-resolution spaces analysis'(MRA) \cite{EPJ} in wavelet analysis, mentioned in Section \eqref{SEC2}. 
Essentially, noting the fact that the idea of multi-resolution spaces was derived from a multi-resolution approximation approach inspired by Mallat and Meyer \cite {E3}.
Roughly and mathematically speaking that the interpolant can be given by a modulated $sinc$ function in the classical Shannon sampling theorem in the form of a scaling function in the multi-resolution spaces.

In \cite{I8,prev}, the reconstruction of function $f(t)$ over the decomposition of  multiresolution spaces \(\{ V_J(\chi),J\in \mathbb{Z}\}\) is given as follows

\begin{eqnarray*}
    f(t) = \frac{1}{2^J} \sum_{\varrho \in \mathbb{Z}} f\left( \frac{\varrho}{2^J} \right) \chi(2^J t - \varrho),
\end{eqnarray*}
and 
\begin{eqnarray*}
    \hat{\chi}(\varepsilon)=\sum_{\rho\in \mathbb{Z}}\frac{\hat{\chi_1}(\varepsilon)}{\hat{\chi_1}(\varepsilon+2\rho \pi)},
\end{eqnarray*}
where $\chi$ and $\hat{\chi}(\varepsilon)$ represent the scaling function (interpolant) and Fourier transform of $\chi(t)$ respectively. Again, noting the fact that 
 in the summation the term $\hat{\chi_1}(\varepsilon)$ is the  Fourier transform of the function $\chi_1(t)$ with specific properties given as follows;
 
\begin{equation} \label{eq:chi_function}
    \chi_1(\varrho)= \begin{cases} 
      0 & \text{if } \varrho= 0, \\
      1 & \text{if } \varrho=\pm 1, \pm 2, \pm 3, \dots
\end{cases}.
\end{equation}

For more detail, the readers are suggested to go through the research paper \cite{ S1}. Consequently, the sampling theorem in terms of wavelets can be  written as follows
 \begin{eqnarray*}
    f(t) = \frac{1}{2^J} \sum_{\varrho \in \mathbb{Z}} f\left( \frac{\varrho}{2^J} \right) \chi_1(2^J t - \varrho),\,\,\,\forall f\in V_J(\chi_1)
\end{eqnarray*}
Moreover, let us denote  $\Tilde{f}(t)$ as a reconstructed function that comes from $V_{J+1}(\chi)$ (not from $V_{J}(\chi)$ ), defined as follows
 \begin{eqnarray*}
     \Tilde{f}(t)=\frac{1}{2^J} \sum_{\varrho \in \mathbb{Z}} f\left( \frac{\varrho}{2^J} \right) \chi(2^J t - \varrho).
 \end{eqnarray*}

\par
In \cite{o1}, authors classify the compactly supported orthonormal scaling functions (COSF) that satisfy equation \eqref{eq:chi_function}  and  establishes a result as; if \( \chi_1(t) \) is a COSF then \( \chi_1(t) \) is the Haar scaling function and converse of this statement is also true. Generally, the Haar wavelet is associated with an indicator function as $\varkappa_1= {\chi_1}_{[0,1]}(t) $ (say). Essentially, the generalization of a Haar function \( {\chi_1}_{[0,1)}(t) \) in a limiting sense can be considered as a collection of a COSF with simple forms and exponential decrement.  Hence, the new form of the reconstruction of a function can be reframed as follows; 

\begin{eqnarray}\label{sampp}
    f(t) = \frac{1}{2^J} \sum_{\varrho \in \mathbb{Z}} f\left( \frac{\varrho}{2^J} \right) \varkappa_1(2^J t - \varrho)
\end{eqnarray}
\par
Inspired by these articles \cite{s5,s6,s3,s2,s1,s4}, Bardaro with his colleague established a Durrmeyer-sampling type operators given by the following expression
\begin{eqnarray*}
     \label{E3}
(\mathfrak{K}_w^{\chi_1,\chi_2}f)(t)=\sum_{\varrho=-\infty}^{+\infty}\left (w\int_{\mathbb{R}} \chi_2(wu-\varrho)du\right)\chi_1(wt-\varrho), t\in \mathbb{R},
\end{eqnarray*}
where the meaning of the parameters $\chi_1,\, \&\,\chi_2$ ar given in section \eqref{SEC2}
Moreover, as discussed in \cite{MAIN} the suitable reconstruction of a function $f(t)$ in the form of wavelet analysis is derived with the help of  equation $\eqref{W1}$ is given as;
$$ f(t)=\digamma_{\upsilon}^{-1} \int_{\mathbb{R}} \int_{\mathbb{R}} (\mathcal{W}_{\Upsilon} f)(a_1,a_2)\Upsilon_{a_1,a_2}(t) \lvert a_1 \rvert ^{2\varrho-3} \, da_1 \, da_2$$
$\text{where}, 0<\digamma_{\Upsilon}^{-1} =\int_{\mathbb{R}}\frac{\widehat{\Upsilon}(w)}{|w|} dw<\infty,\, t \in \mathbb{R},\, \Upsilon_{a_1,a_2}=a_1^{\varrho} \Upsilon\left(\frac{t-a_2}{a_1}\right)$ and \,$(\mathcal{W}_{\Upsilon} f)(a_1,a_2)$ denotes the wavelets transform. \par
\vskip0.15in

  Later on, Hirschman \cite{I12} derives a convolution property associated with wavelet transform in terms of a basic function $D(t_1, t_2, t_3)$, which is given as follows;
\begin{eqnarray*}
    (\mathcal{W}_{\Upsilon} \,(g^* \Xi g)\,)(a_1,a_2)  &= & \int_{\mathbb{R}} \overline{\chi_2(2^{j}u-\varrho)}\,d u \int_{\mathbb{R}} \int_{\mathbb{R}} D(u,t_2,t_3) f(t_3) dt_2 dt_3\\ 
    &= &\int_{\mathbb{R}} (g^* \Xi g)(u) \overline{\chi_2(2^{j}u-\varrho)}\,d u, \\
\end{eqnarray*}
where the terms like  $g\,,\, g^* \&\, D(u,t_2,t_3)$ are defined in section $2.$
  \par

  Thus, the `wavelet-based filtering operators' using \eqref{sampp} is given as follows:
\begin{equation}\label{E4}
( \mathfrak{W}_{\varrho}^{\varkappa_1,\chi_2} h )(t) = \sum_{\varrho=-\infty}^{+\infty} \varkappa_1(2^{j} t - \varrho) \int_{\mathbb{R}} h(u) \chi_2(2^{j} u - \varrho) \, du
\end{equation}

where $ a_1\in \mathbb{R}_+, a_2\in \mathbb{R}, \varrho\geq 0, \Upsilon, \, \theta, \upsilon\, $and $ h\in L^1(\mathbb{R})$ in terms of the basic function $D(t_1, t_2, t_3)$  is defined in \eqref{h(t)} as $ \chi_2(t)=(f'\, \Xi\,\, \chi_2')(t)=\int_{\mathbb{R}} D(t, t_2, t_3)  \chi_2(t_2)\,dt_2\, dt_3$. We will study these terms in detail in section \eqref{SEC2}.\par
The main focus of the article entails the construction of the `wavelet-based filtering operators' and using the basic definitions and terminologies, the core of the article in the forms of theorems and prepositions is derived in \eqref{E4}. \eqref{SEC2} and \eqref{Sec3}.  Eventually,  some examples are illustrated, giving  better clarity to understand the theory as mentioned in \eqref{EXA}.

\vskip0.15in

 \section{Preliminaries}\label{SEC2}
The notation \( (\cdot, \cdot) \) and \( \| \cdot \| \) represent the inner product and the norm in \( L^2(\mathbb{R}) \), respectively, where \( L^2(\mathbb{R}) \) denotes the space of all square-integrable functions over the real numbers \( \mathbb{R} \). Furthermore, let us suppose that \( \hat{f}(\omega) \) denotes the Fourier transform of a function \( f(t)\) which is square integrable functions, written as follows;
\[
\hat{f}(\omega) = \int_{-\infty}^{\infty} f(t) e^{-i \omega t} \, dt.
\]

Let us denote the space of square-summable sequences by \( \ell^2(\mathbb{Z}) \) over the set of integers $\mathbb{Z}$. In the context of wavelet approximation, we define the scaling function \( \chi(t) \), which is summable sequence over $\mathbb{Z}$ and throughout the paper, we assume that it is orthogonal.

In order to define the multiresolution approximation (MRA), for $j\in \mathbb{Z}$, we have subspaces \( V_j(\cdot) \) such that \( V_j(\cdot) \subseteq V_{j+1}(\cdot) \) with a wavelet function $\Upsilon(t) \in  L^2(\mathbb{R})$. Further, let  \( \frac{1}{2}{h_1}_{\varrho} \) and \( \frac{1}{2}{h_2}_{\varrho} = (-1)^{\varrho} \frac{1}{2} p_{1-{\varrho}} \)be the impulse responses of low-pass filter \( L(\omega) \) and \( B(\omega) \) respectively, with respect  to the scaling function.\\
Now, we study some important properties on the basis of scaling functions and impulse responses as follows;\par
\begin{itemize}

  \item Dilation Equation:
\begin{eqnarray*}
    \chi_1(t)=\sum_{\varrho \in \mathbb{R}} {h_1}_{\varrho} \chi_1(2t-\varrho)
\end{eqnarray*}
$\& $
\begin{eqnarray*}
    \Tilde{\chi_1}(\omega)=\sum_{\varrho \in \mathbb{Z}} {h_1}_\varrho (\omega2^\varrho)
\end{eqnarray*}
\end{itemize}

\begin{itemize}

    \item  Orthogonality:
\begin{eqnarray}\label{1}
    |{h_1}_\varrho(\omega)|^2+ |{h_1}_\varrho(\omega+\pi)|^2=1,\, \omega\in \mathbb{R} 
\end{eqnarray}
\end{itemize}
\begin{itemize}

    \item  Normality:
\begin{eqnarray}\label{2}
    \sum_{\varrho\in\mathbb{Z}}{h_1}_\varrho=2
\end{eqnarray}
\end{itemize}

The wavelet $ \chi_2(t)$ can be constructed by
\begin{eqnarray*}
    \Upsilon(t)=\sum_{\varrho \in \mathbb{Z}} {h_2}_\varrho \chi_2(2t-\varrho)
\end{eqnarray*}
In \cite{o1,o2}, the smoothness of the  scaling functions  \eqref{1} and the wavelets \eqref{2}  depends on the impulse response \( {h_1}_\varrho \) corresponding to  low-pass filter $ L(\omega)$ satisfying the following conditions ;

\begin{eqnarray}\label{3}
    \left\lvert\frac{d^\varrho L(\omega)}{d \omega^\varrho}\right\rvert_{\omega=\pi}=1,\,\, 1\leq \varrho \leq I-1
\end{eqnarray}
for certain integer $I\geq 1$. If ${h_1}_{\varrho}$ satisfies \eqref{3}, then the 
constructed functions have $q_{th}$ (say) order smoothness, where $q \propto$ I.
\begin{itemize}
    
\item The modulus of continuity for a function 
$f$ on  $[a_1,a_2]$ is defined as:
\begin{eqnarray*}
      \omega(f, \delta) = \sup_{\substack{x, y \in [a_1, a_2] \\ |x - y| \leq \delta}} |f(x) - f(y)|. 
\end{eqnarray*}
There are various theorems and propositions in the form of estimation and rate of convergence are derived with the help of madulus of continuity $\omega(f, \delta)$ as mentioned in the articles\cite{R3,V1,V2,R1,R2,DS}.\\
\end{itemize}

Let us assume that wavelets and impulse responses \(\chi_2(t)\), \(\Upsilon(t)\), \({h_1}_\varrho\) and \({h_2}_\varrho\), satisfying $(2.1)$ and $(2.2)$. Suppose that the sample values \(t[\varrho]\) for \(J \geq 0\), with \(t[\varrho] = f(\varrho/2^J)\). Again, let us reframe  \(c_{J,\varrho}\) and \(b_{j,\varrho}\) in a recursive way as defined in \eqref{7} and \eqref{8}, respectively in the following form with the condition on \(b_{J,\varrho}\), for \(j\leq J.\), \(\varrho \in \mathbb{Z}\), and \(c_{J,\varrho}\), \(k \in \mathbb{Z}\);
\begin{eqnarray}\label{4}
     c_{j-1,\varrho}=\frac{1}{2} \sum_{\varrho\in \mathbb{Z}}{h_1}_{\varrho-2\varrho} c_{j,n}
\end{eqnarray}
\begin{eqnarray}\label{5}
     b_{j-1,\varrho}=\frac{1}{2} \sum_{\varrho\in \mathbb{Z}} {h_2}_{\varrho-2\varrho} b_{j,n}
\end{eqnarray}
for $j<J.$  For our convenience, we denote a discrete signal by \(y_d[\varrho], \varrho \in \mathbb{Z}\), as initial sequence as defined in equations \eqref{4} and \eqref{5}, instead of using the sequence \(c_{J,\varrho}, \varrho \in \mathbb{Z}\).  Equivalently, the decomposition \(b_{j,\varrho}, j \leq J - 1, \varrho \in \mathbb{Z}\) refers to \(\text{DWT}\{y_d[\varrho]; j,\varrho \}\), where DWT stands for discrete wavelet transform. If the initial sequence \(c_{J,\varrho}, \varrho \in \mathbb{Z}\), in \eqref{4} and \eqref{5} corresponds to the discrete-time signal \(t[\varrho], \varrho \in \mathbb{Z}\), sampled from \(f(t)\), then by Mallat algorithm, WST coefficients of \(f(t)\) is equivalent to \(\text{DWT}\{ t'[\varrho],; j, \varrho, j \leq J - 1, \varrho \in \mathbb{Z}\)\}. During the discussion of the Shensa algorithm, we use \( t'[\varrho]\), a prefiltered version of \( t[\varrho]\) (as in the Mallat algorithm), as the initial sequence in \eqref{4} and \eqref{5} to compute \(b_{j,\varrho}\), for \(j < J, \varrho \in \mathbb{Z}\), we have 

\begin{eqnarray*}
    t'[\varrho] = \frac{1}{2} \sum_{s \in \mathbb{Z}} t[s] p[\varrho - s]
\end{eqnarray*}

\begin{eqnarray*}
    p[\varrho]=2^{-\frac{J}{2}}\int_{\mathbb{R}} \chi_1(t)p(t-\varrho) dt
\end{eqnarray*}
Similarly, the WST coefficients of \( f(t) \) namely \( b_{J,\varrho}, \, j \leq J - 1, \, \varrho \in \mathbb{Z} \) corresponding to the interpolant of \( t[\varrho] \) are calculated obtained using the Shensa algorithm and are denoted by \( b_{j,\varrho}^* \). Further, assume that if we denote the signal \( \Tilde{f}(t) \) produced by interpolant  \( t[\varrho] \), then
\begin{eqnarray*}
    \Tilde{f}(t) = \frac{1}{2^J} \sum_{\varrho \in \mathbb{Z}} t^{'}[\varrho] \varkappa_1(2^J t - \varrho)
\end{eqnarray*}

\begin{eqnarray}\label{samp}
    \Tilde{f}(t) = \frac{1}{2^J} \sum_{\varrho \in \mathbb{Z}} f\left( \frac{\varrho}{2^J} \right) \varkappa_1(2^J t - \varrho)
\end{eqnarray}

In \cite{w10}, Hirschman develops a novel technique using basic function $D(t_1, t_2, t_3)$ defined in \eqref{W1} and the convolution associated with wavelet transform, defined as follows;

\begin{eqnarray*}
     (\mathcal{W}_{\Upsilon} D(., t_2, t_3))(a_1,a_2)&=&\int_{\mathbb{R}} D(t_1, t_2, t_3) \overline{\chi_2(2^{j} t - k})\,d t_1\\
     &=& \overline{\Upsilon}(2^{j} t_3 - \varrho) \overline{\chi_2}(2^{j} t_2 - \varrho)
\end{eqnarray*}
where, $$ D(t_1, t_2, t_3)=\digamma_{\upsilon}^{-1}\int_{\mathbb{R}} \int_{\mathbb{R}} \overline{\Upsilon}(2^{j} t_3 - \varrho) \overline{\chi_2}(2^{j} t_2 - \varrho) \upsilon(2^{j} t_1 - \varrho)\lvert a_1 \rvert ^{2\varrho-3} \, da_1 \, da_2 $$
Let the translation $\tau_{t_1}$ is defined in terms of wavelet structure which is given as follows
\begin{eqnarray*}
    (\tau_{t_1} f)(t_2) &=& g^*(t_1,t_2)= \int_{\mathbb{R}} D(t_1, t_2, t_3) f( t_3)\, dt_3\\
    &=&\digamma_{\upsilon}^{-1}\int_{\mathbb{R}}\int_{\mathbb{R}} \int_{\mathbb{R}} \overline{\Upsilon}(2^{j} t_3 - \varrho) \overline{\theta}(2^{j} t_2 - \varrho) \upsilon(2^{j} t_1 - \varrho) f( t_3)\lvert a_1 \rvert ^{2\varrho-3} \, da_1 \, da_2 \, dt_3 
\end{eqnarray*}
Furthermore, let us denote the convolution by  $\Xi$ between $g^*(t_1,t_2)$ and $ g(2^{j} t_2 - \varrho)$ defined in a following way;
\begin{eqnarray*}
(g^* \Xi\,\, g)(t)&= &\int_{\mathbb{R}} g*(t,t_2)g(t_2)\, dt_2\\
&= &  \int_{\mathbb{R}}\int_{\mathbb{R}} D(t, t_2, t_3)   g( t_2 ) f(t_3 )\,dt_2\, dt_3\\
&= &\digamma_{\upsilon}^{-1}\int_{\mathbb{R}}\int_{\mathbb{R}} \int_{\mathbb{R}} \overline{\Upsilon}(2^{j} t_3 - \varrho) g( t_2) \overline{\theta}(2^{j} t_2 - \varrho)\upsilon(2^{j} t_1 - \varrho) f( t_3 )\\
\end{eqnarray*}
\begin{eqnarray}\label{eq2.4}
&\times&\lvert a_1 \rvert ^{2\varrho-3} \, da_1 \, da_2 \,dt_2 dt_3 
\end{eqnarray}

\vskip0.15in
\textbf{Proposition:1} The sampling theorem stated as in \eqref{samp} is valid if $f(t)$ is a cardinal scaling function and vice versa.\\ 
\vskip0.15in
\textbf{Note 1} Specifically, if \(a_1 = 2^{-j}\), \(a_2 = k 2^{-J}\), where \(j < J\) and \(k \in \mathbb{Z}\), for some \(J > 0\), then $da_1=-2^{-j}\log2dj\,\&\,\,\,da_2=-2^{j}k\log2 dj$ then from equation\eqref{eq2.4} we have
\begin{eqnarray*}
\end{eqnarray*}
\begin{eqnarray*}
\chi_2(t)=(g^{*'} \Xi\,\, g')(t)
&= &  \int_{\mathbb{R}}\int_{\mathbb{R}} D'(t, t_2, t_3)   g'
( t_2 ) f'(t_3 )\,dt_2\, dt_3\\
&= &\lvert 2 \rvert ^{(3J-2J\varrho)+1}\,\varrho\log2\digamma_{\upsilon}^{-1}\int_{\mathbb{R}}\int_{\mathbb{R}} \int_{\mathbb{R}} \overline{\Upsilon}(2^{j} t_3 - \varrho) \overline{\chi_2}(2^{j} t_2 - \varrho) \upsilon(2^{j} t_1 - \varrho)\\
\end{eqnarray*}
\begin{eqnarray}\label{h(t)}
    &&\times dj^2dt_2 dt_3 
\end{eqnarray}

\textbf{Note 2} The reconstruction of $f(t)$ in terms of wavelets $\Upsilon_{(2^{j} t - \varrho)}$ is given by $f(t)=\chi_{J,\varrho}=\Upsilon(2^{j} t - \varrho)$ with
\begin{eqnarray*}
    f(t)=\sum_{j,\varrho}b_{J,\varrho}\Upsilon_{J,\varrho}
\end{eqnarray*}
\begin{eqnarray}\label{7}
    b_{j,\varrho}=<f, \Upsilon_{J,\varrho}>
\end{eqnarray}
\textbf{Note 3} The notation $f_j(t)$ denotes the orthogonal projection in $V_j(\chi_1)$ of $f(t)$, then
\begin{eqnarray*}\label{8}
    f_j(t)=\sum_{k}c_{J,\varrho}\chi_{J,\varrho}(t)=\sum_{j<J}\sum_{\varrho}b_{J,\varrho}\Upsilon_{J,\varrho}(t),
\end{eqnarray*}
where $\chi_{J,\varrho}=\chi_1(2^{j} t - \varrho)$ and
\begin{eqnarray}
    c_{j,\varrho}=<f, \chi_{J,\varrho}>
\end{eqnarray}
\begin{theorem}\label{TH0}
    Suppose that $\upsilon \in L^p(\mathbb{R}), \theta \in L^q(\mathbb{R}) \,\, \&  \,\, \Upsilon \in L^1(\mathbb{R})$ with the condition such that $\frac{1}{p}+\frac{1}{q}=1+\varrho\,,\, 0<\varrho<1\&\,\, p,q>1$
\begin{eqnarray*}
     \int_{\mathbb{R}}|D(t_1, t_2, z)|dz \leq \digamma_{\upsilon}^{-1} K(p,\varrho)  |t_1-t_2|^{-\varrho}\lVert \theta \rVert_q \lVert \Upsilon \rVert_1 \lVert\upsilon\rVert_p , 
\end{eqnarray*}
where $K(p,\varrho)$ is a constant.\\
\end{theorem} \label{TH1}
\begin{theorem}\label{TH1}
    If \((1 + |t|\varrho)\upsilon(t) \in L^1(\mathbb{R})\) and \(\Upsilon \in L^1(\mathbb{R})\) such that \(\varrho\geq 1\). Further, assume that \((1 + |t|\varrho)\theta(t) \in L^1(\mathbb{R})\). Then
\begin{eqnarray*}
     \int_{\mathbb{R}}|D(t_1, t_2, z)|dz \leq 2^{\varrho -1}\digamma_{\upsilon}^{-1} G(p,\varrho)  |t_1-t_2|^{-\varrho}\left[\lVert \upsilon(t) t^{\varrho-1} \rVert_1\lVert \theta(t) \rVert_1 +\lVert \theta(t) t^{\varrho-1} \rVert_1\lVert \lVert \upsilon(w) \rVert_1\right]\Upsilon(w)\rVert_1  , 
\end{eqnarray*}
\end{theorem} 
For more details and proofs, the keen readers are suggested to go through the article \cite{w10}.

\section{Main Results} \label{Sec3}
Now, we will prove the main result based on approximation with respect to the operators defined in $\eqref{E4}$
\subsection{Approximation by $\mathfrak{W}_{\varrho}^{\varkappa_1,\chi_2}h$ operators
}
\begin{theorem}
    Suppose $h \in \mathbb{C}^N(\mathbb{R})$ and $\varrho \in \mathbb{Z}$. Further, suppose that $\varkappa_1$ is a bounded Haar scaling function with compact support in $[-a, a]$ for some $a > 0$, such that
\[
\sum_{\varrho=-\infty}^{\infty} \varkappa_1(t - \varrho) = 1, \quad \text{for all } t \in \mathbb{R}.
\]
Moreover, suppose $\varkappa_1 \geq 0$. Then 
\label{T1}
\begin{eqnarray}
    \lvert ( \mathfrak{W}_{\varrho}^{\varkappa_1,\chi_2}h)(t)-h(t)\rvert \leq \sum_{i=1}^{\varrho} \frac{\lvert h^i(x)\rvert}{i!}\frac{a^i}{2^{i\varrho}}+ \frac{a^\mathbb{N}}{2^{\mathbb{N}\varrho}}\omega_1\left( h^{\mathbb{N}}, \frac{a}{2^{j}}\right), 
\end{eqnarray}
where the operators $\mathfrak{W}_{\varrho}^{\varkappa_1,\chi_2}h$ is defined in $\eqref{E4}.$
\proof
Given \( h \in L^p(\mathbb{R}) \), the operators introduced in \eqref{E4} provide the following outcomes:

\begin{eqnarray*}
    \lvert ( \mathfrak{W}_{\varrho}^{\varkappa_1,\chi_2} h)(t) - h(t) \rvert 
    &=& \left\lvert 2^{\frac{j}{2}} \sum_{\varrho =-\infty}^{+\infty} \left[ \int_{\mathbb{R}} h(u) \chi_2(2^{j} u - \varrho) \, du - 2^{-\frac{j}{2}} h(t)\right] \varkappa_1(2^{j} t - \varrho) \right\rvert \\
    &=& \left\lvert 2^{\frac{j}{2}} \sum_{\varrho =-\infty}^{+\infty} \left[ \int_{\mathbb{R}} h\left(\frac{v}{2^{j}}\right) \chi_2(v - k) \, dv - 2^{-\frac{j}{2}} h(t)\right] \varkappa_1(2^{j} t - k)\right\rvert \\
    &=& \left\lvert 2^{\frac{j}{2}} \sum_{\varrho =-\infty}^{+\infty} \left[ \int_{\mathbb{R}} h\left(\frac{v}{2^{j}}\right) \,dv - h(t)\right]\chi_2(v - k) \varkappa_1(2^{j} t - k)\right\rvert \\
\end{eqnarray*}

As $sup \varkappa_1 \in [-a,a]$,\,then \, $\forall t\,\in \mathbb{R}$, we have $2^{j}{(-a+j)}\leq t \leq 2^{j}(a+j).$ Thus we have
\begin{eqnarray*}
    \lvert ( \mathfrak{W}_{\varrho}^{\varkappa_1,\chi_2} h)(t) - h(t) \rvert 
    &=& \left\lvert 2^{\frac{\varrho}{2}}\Delta \right\rvert,
\end{eqnarray*}
where $\Delta=\sum_{\varrho =-\infty}^{+\infty} \left[ \int_{-a+j}^{a+j} h\left(\frac{v}{2^{j}}\right) \,dv - h(t)\right]\chi_2(v - \varrho) \varkappa_1(2^{j} t - \varrho) $.\\
\begin{eqnarray*}
    h\left(\frac{v}{2^{j}}\right)  - h(t)&=& \sum_{i=1}^{\mathbb{N}}\frac{f^i(t)}{i!}\left(\frac{v}{2^{j}}--t\right)^i\\
    & & +\int_{t}^{\frac{v}{2^{j}}}\left ( h^{\mathbb{N}}\left(\frac{s}{2^{j}}\right) \,dv - h^{\mathbb{N}}(t)\right) \frac{\left(\frac{s}{2^{j}}-t\right)^{\mathbb{N}-1}}{(\mathbb{N}-1)!}ds\\
 \implies   \left[h\left(\frac{v}{2^{j}}\right) \,dv - h(t)\right]\chi_2(v - k) \chi_1(2^{j} t - k)&=&\sum_{i=1}^{\mathbb{N}}\frac{f^i(t)}{i!}\left(\frac{v}{2^{j}}-t\right)^i \chi_2(v - k) \varkappa_1(2^{j} t - k)\\
    &+& \chi_2(v - k) \varkappa_1(2^{j} t - k) \times\\
    && \int_{t}^{\frac{v}{2^{j}}}\left ( h^{\mathbb{N}}\left(\frac{s}{2^{j}}\right) \,dv - h^{\mathbb{N}}(t)\right) \frac{\left(\frac{s}{2^{j}}-t\right)^{\mathbb{N}-1}}{(\mathbb{N}-1)!}ds.\\
\end{eqnarray*}
Clearly, $\left\lvert \frac{2}{\varrho}-t\right\rvert < \frac{2}{\varrho-1}$, we have
\begin{eqnarray*}
    \Delta &=& \sum_{i=1}^{\mathbb{N}}\frac{f^i(t)}{i!}\left(\frac{v}{2^{j}}-t\right)^i \chi_2(v - \varrho) \varkappa_1(2^{j} t - \varrho)+\mathbb{E_R},\\
\end{eqnarray*}
where 
\begin{eqnarray*}
    \mathbb{E_R} &=& \chi_2(v - \varrho) \varkappa_1(2^{j} t - \varrho) \int_{t}^{\frac{v}{2^{j}}}\left ( h^{\mathbb{N}}\left(\frac{s}{2^{j}}\right) \,dv - h^{\mathbb{N}}(t)\right) \frac{\left(\frac{s}{2^{j}}-t\right)^{\mathbb{N}-1}}{(\mathbb{N}-1)!}ds.
\end{eqnarray*}
Conclusively, $|\Delta|\leq \sum_{i=1}^{\mathbb{N}}\frac{f^i(t)}{i!}\frac{a^i}{2^{i(\varrho-1)}} +|\mathbb{E_R}|, $
where $|\mathbb{E_R}|$ is estimated as follows
\begin{eqnarray*}
    |\mathbb{E_R}|&\leq& \varkappa_1(2^{j} t - \varrho)\int_{-a+j}^{a+j} \chi_2(v - \varrho)  \xi (v) dv\\
\end{eqnarray*} 
with the condition $t\leq \frac{v}{2^{j}}$, we have
\begin{eqnarray*}
     \xi(v) &\leq & \left \lvert\int_{t}^{\frac{v}{2^{j}}}\left ( h^{\mathbb{N}}\left(s\right) - h^{\mathbb{N}}(t)\right) \frac{\left(\frac{v}{2^{j}}-s\right)^{\mathbb{N}-1}}{(\mathbb{N}-1)!}ds\right \rvert\\
     &\leq &\int_{t}^{\frac{v}{2^{j}}} \omega_1(f^{\mathbb{N}},|s-t|)\frac{\left(\frac{v}{2^{j}}-s\right)^{\mathbb{N}-1}}{(\mathbb{N}-1)!}ds\\
      &\leq & \omega_1(f^{\mathbb{N}},\frac{a}{2^{\varrho-1}})\frac{\left(\frac{v}{2^{j}}-t\right)^{\mathbb{N}}}{(\mathbb{N})!}\\
       &\leq & \omega_1(f^{\mathbb{N}},\frac{a}{2^{\varrho-1}})\frac{\left(\frac{v}{2^{j}}-t\right)^{\mathbb{N}}}{2^{\varrho-1}(\mathbb{N})!}\\
\end{eqnarray*}
Now, if $t\geq \frac{v}{2^{j}}$ we have
\begin{eqnarray*}
  \xi(v)  &\leq & \left \lvert\int_{\frac{v}{2^{j}}}^t\left ( h^{\mathbb{N}}\left(s\right) - h^{\mathbb{N}}(t)\right) \frac{\left(\frac{v}{2^{j}}-s\right)^{\mathbb{N}-1}}{(\mathbb{N}-1)!}ds\right \rvert\\
   &\leq & \omega_1(f^{\mathbb{N}},\frac{a}{2^{\varrho-1}})\frac{\left(t-\frac{v}{2^{j}}\right)^{\mathbb{N}}}{(\mathbb{N})!}\\
   &\leq & \omega_1(f^{\mathbb{N}},\frac{a}{2^{\varrho-1}})\frac{\left(t-\frac{v}{2^{j}}\right)^{\mathbb{N}}}{2^{\varrho-1}(\mathbb{N})!}\\
\end{eqnarray*}
The estimation of $\xi(v)$ in both the cases are same. Therefore, let us assume this $W$
\begin{eqnarray*}
    |\mathbb{E_R}|&\leq& \left(\varkappa_1(2^{j} t - \varrho)\int_{-a+j}^{a+j} \chi_2(v - \varrho)  \xi (v) dv\right)W=:W\\
\end{eqnarray*} 

\begin{eqnarray*}
   |\mathbb{E_R}|&\leq&\omega_1(f^{\mathbb{N}},\frac{a}{2^{\varrho-1}})\frac{\left(t-\frac{v}{2^{j}}\right)^{\mathbb{N}}}{2^{\varrho-1}(\mathbb{N})!}\\
\end{eqnarray*} 
\begin{eqnarray*}
    \Delta &=& \sum_{i=1}^{\mathbb{N}}\frac{f^i(t)}{i!}\left(\frac{v}{2^{j}}-t\right)^i \chi_2(v - \varrho) \varkappa_1(2^{j} t - \varrho)+\omega_1(f^{\mathbb{N}},\frac{a}{2^{\varrho-1}})\frac{\left(t-\frac{v}{2^{j}}\right)^{\mathbb{N}}}{2^{\varrho-1}(\mathbb{N})!}=\Theta,\\
\end{eqnarray*}
Thus, we have, 
$$|\Delta |2^{\frac{-\varrho}{2}}\leq| \Theta| 2^{\frac{-\varrho}{2}}$$
\begin{eqnarray*}
\left[ \int_{\mathbb{R}} h(u) \chi_2(2^{j} u - \varrho) \, du - 2^{-\frac{k}{2}} h(t)\right]\leq \Theta 2^{\frac{-\varrho}{2}}
\end{eqnarray*}
More preciously, 
\begin{eqnarray*}
    \lvert ( \mathfrak{W}_{\varrho}^{\varkappa_1,\chi_2} h)(t) - h(t) \rvert 
    &=& \left\lvert 2^{\frac{j}{2}} \sum_{\varrho =-\infty}^{+\infty} \left[ \int_{\mathbb{R}} h(u) \chi_2(2^{j} u - \varrho) \, du - 2^{-\frac{J}{2}} h(t)\right] \varkappa_1(2^{j} t - \varrho) \right\rvert\\
     &\leq& \Theta \sum_{\varrho =-\infty}^{+\infty}  \varkappa_1(2^{j} t - \varrho)= \Theta . \\
\end{eqnarray*}
\end{theorem}
\textbf{Note 4.} Let us assume that the function \( h^{\mathbb{N}} \) is a bounded continuous function over \( \mathbb{R} \).  
Then, the modulus of continuity \( \omega\left(h^{\mathbb{N}}, \frac{a}{2^{\varrho-1}}\right) \) is finite, and  
\[
\mathfrak{W}_{\varrho}^{\varkappa_1, \chi_2} (h)(t) \to h(t) \quad \text{as} \quad \varrho \to +\infty, \quad \forall t \in \mathbb{R}.
\]

\subsection{Estimation of a function $h\in L^{r_1}(\mathbb{R})\,\cap\in \L^{r_2}(\mathbb{R})$ with respect to the operators $\mathfrak{W}_{\varrho}^{\varkappa_1,\chi_2} h$} 

\begin{theorem}\textbf{(a)} Suppose that $\upsilon \in L^p(\mathbb{R}), \theta\in L^q(\mathbb{R}) \,\, \&  \,\, \Upsilon \in L^1(\mathbb{R})$ with the condition such that $\frac{1}{p}+\frac{1}{q}=1+\varrho\,,\, 0<\varrho<1\&\,\, p,q>1$. Furthermore,  $h(t)=f\,\Xi \chi_2,\,\,,f\in L^{r_1}(\mathbb{R})\,\,\chi_2\in \L^{r_2}(\mathbb{R}) , $ such that $ \Upsilon \in L^q(\mathbb{R}), \, \chi_2 \in L^1(\mathbb{R})$ with $\frac{1}{p}+\frac{1}{q}=1+\nu,p,q>1\, \nu \in (0,1)\,\, \frac{1}{r_1}+\frac{1}{r_2}+\nu=2$. Then 
\begin{eqnarray*}
     \lvert ( \mathfrak{W}_{\varrho}^{\varkappa_1,\chi_2} h)(t) - h(t) \rvert&\leq &\lvert 2 \rvert ^{(3J-2J\varrho)+1}\,\varrho\log2\digamma_{\upsilon}^{-1} K(p,\varrho) \lVert \theta \rVert_q \lVert \Upsilon \rVert_1 \lVert\upsilon\rVert_p K^{'}(p,\varrho)\\
&&\times\lVert \theta \rVert_q \lVert \Upsilon \rVert_1 \lVert\upsilon\rVert_p f(t_3) \rVert_r \lVert \chi_2 \rVert_{r'}.\\
\end{eqnarray*}
where $K$ is a real constant.
\proof
Firstly, we will derive the result for the $0<\varrho < 1$. As given $h(t)=f'\,\Xi \chi_2',\,\,,f'\in L^{r_1}(\mathbb{R})\,\,\chi_2\in \L^{r_2}(\mathbb{R}) $, the operators defined \eqref{E4}, we have
    \begin{eqnarray*}
    \lvert ( \mathfrak{W}_{\varrho}^{\varkappa_1,\chi_2} h)(t) - h(t) \rvert 
    &=& \left\lvert 2^{\frac{J}{2}} \sum_{\varrho =-\infty}^{+\infty} \left[ \int_{\mathbb{R}} h(u) \chi_2(2^{j} u - \varrho) \, du - 2^{-\frac{J}{2}} h(t)\right] \chi_1(2^{j} t - \varrho) \right\rvert \\
      &=& \left\lvert 2^{\frac{J}{2}} \sum_{\varrho =-\infty}^{+\infty} \left[ \int_{\mathbb{R}} (f'\Xi \chi_2')(t,u)\chi_2(2^{j} u - \varrho) \, du - 2^{-\frac{J}{2}} h(t)\right]\right.\\
      &&\times \left. \chi_1(2^{j} t - \varrho) \right\rvert \\
      &\leq & \left\lvert 2^{\frac{J}{2}} \sum_{\varrho =-\infty}^{+\infty} \left[ \left( \int_{\mathbb{R}}\int_{\mathbb{R}} D(t, u, t_3)  \chi_2(2^{j} u - \varrho)f(t_3)\,du\, dt_3\right)  \, du - 2^{-\frac{J}{2}} h(t)\right]\right.\\
      &&\times \left. \chi_1(2^{j} t - \varrho) \right\rvert \\
    \end{eqnarray*}
    From equation \eqref{TH0}, we obtain 
    \begin{eqnarray*}
 \lvert ( \mathfrak{W}_{\varrho}^{\varkappa_1,\chi_2} h)(t) - h(t) \rvert &\leq & \lvert 2 \rvert ^{(3J-2J\varrho)+1}\,\varrho\log2\digamma_{\upsilon}^{-1} K(p,\varrho) \lVert \theta \rVert_q \lVert \Upsilon \rVert_1 \lVert\upsilon\rVert_p   \\
&\times&  \left|\sum_{\varrho =-\infty}^{+\infty} \left[ \chi_1(2^{j} t - k)\int_{\mathbb{R}} |t-u|^{-\varrho}\chi_2(2^{j} u - \varrho) \, du\int_{\infty} f(t_3) \, dt_3 \right] - 2^{-\frac{k}{2}} h(t) \right|.
\end{eqnarray*}
Now, from Hardy–Littlewood–Sobolev inequality \cite{In1}, we have
\begin{eqnarray*}
\lvert ( \mathfrak{W}_{\varrho}^{\varkappa_1,\chi_2} h)(t) - h(t) \rvert &\leq &\lvert 2 \rvert ^{(3J-2J\varrho)+1}\,\varrho\log2\digamma_{\upsilon}^{-1} K(p,\varrho) \lVert \theta \rVert_q \lVert \Upsilon \rVert_1 \lVert\upsilon\rVert_p  \\
&& \times \left| \sum_{\varrho =-\infty}^{+\infty} \left[ \int_{\mathbb{R}} \int_{\mathbb{R}} \frac{f(t_3)\chi_2(2^{j} u - k)}{|t-u|^{\varrho}} \, du \, dt_3 \right] \right|\lvert \varkappa_1(2^{j} t - \varrho) \rvert \, (2^{\frac{\varrho}{2}}+1)\\
&\leq & \lvert 2 \rvert ^{(3J-2J\varrho)+1}\,\varrho\log2\digamma_{\upsilon}^{-1} K(p,\varrho) \lVert \theta \rVert_q \lVert \Upsilon \rVert_1 \lVert\upsilon\rVert_p K(p,\varrho) \lVert \Upsilon \rVert_q \lVert \Upsilon \rVert_1 \lVert \upsilon \rVert_p \lVert f(t_3) \rVert_r\\
&&\times\lVert \chi_2 \rVert_{r'} \sum_{\varrho =-\infty}^{+\infty} \lvert \varkappa_1(2^{j} t - \varrho) \rvert \, (2^{\frac{J}{2}} + 1)\\
&\leq &\lvert 2 \rvert ^{(3J-2J\varrho)+1}\,\varrho\log2\digamma_{\upsilon}^{-1} K(p,\varrho) \lVert \theta \rVert_q \lVert \Upsilon \rVert_1 \lVert\upsilon\rVert_p K(p,\varrho) \lVert \theta \rVert_q \lVert \Upsilon \rVert_1 \\
&&\times \lVert\upsilon\rVert_p f(t_3) \rVert_r \lVert \chi_2 \rVert_{r'}\, (2^{\frac{J}{2}} + 1)\\
&\leq & \lvert 2 \rvert ^{(3J-2J\varrho)+1}\,\varrho\log2\digamma_{\upsilon}^{-1} K(p,\varrho) \lVert \theta \rVert_q \lVert \Upsilon \rVert_1 \lVert\upsilon\rVert_p K^{'}(p,\varrho)\\
&&\times\lVert \theta \rVert_q \lVert \Upsilon \rVert_1 \lVert\upsilon\rVert_p f(t_3) \rVert_r \lVert \chi_2 \rVert_{r'}.\\
\end{eqnarray*}
where $K'(p,\varrho)$ is a real constant.
\begin{theorem}\textbf{(b)} 
 Suppose that $(1 + |t|\varrho)\upsilon(t) \in L^1(\mathbb{R}), (1 + |t|\varrho)\theta(t) \in L^1(\mathbb{R}) \,\, \&  \,\, \Upsilon \in L^1(\mathbb{R})$ with the condition such that $\frac{1}{p}+\frac{1}{q}=1+\varrho\,,\, \varrho\geq1\&\,\, p,q>1$. Furthermore,  $h(t)=f\,\Xi \chi_2,\,\,,f\in L^{r_1}(\mathbb{R})\,\,\chi_2\in \L^{r_2}(\mathbb{R}) , $ such that $ \Upsilon \in L^q(\mathbb{R}), \, \chi_2 \in L^1(\mathbb{R})$ with $\frac{1}{p}+\frac{1}{q}=1+\nu,p,q>1\, \nu \in (0,1)\,\, \frac{1}{r_1}+\frac{1}{r_2}+\nu=2$. Then 
\begin{eqnarray*}
     \lvert ( \mathfrak{W}_{\varrho}^{\varkappa_1,\chi_2} h)(t) - h(t) \rvert&\leq & 2^{\varrho -1}\digamma_{\upsilon}^{-1} G(p,\varrho) \left[\lVert \upsilon(t) t^{\varrho-1} \rVert_1\lVert \theta(t) \rVert_1 +\lVert \theta(t) t^{\varrho-1} \rVert_1\lVert \lVert \upsilon(w) \rVert_1\right]\\
  &\times& \Upsilon(w)\rVert_1f(t_3) \rVert_r \lVert \chi_2 \rVert_{r'},\\
\end{eqnarray*}
where $K^{'}(p,\varrho) $ is a real constant in terms of $p\,\&\,\varrho\in \mathbb{R}$.
    \proof Let us calculate the estimate between the operators defined in \eqref{E4} and the function $h(t)$, we have \begin{eqnarray*}
   \lvert ( \mathfrak{W}_{\varrho}^{\varkappa_1,\chi_2} h)(t) - h(t) \rvert 
      &\leq & \left\lvert 2^{\frac{J}{2}} \sum_{\varrho =-\infty}^{+\infty} \left[ \left( \int_{\mathbb{R}}\int_{\mathbb{R}} D(t, u, t_3)  \chi_2(2^{j} u - \varrho)f(t_3)\,du\, dt_3\right)  \, du - 2^{-\frac{J}{2}} h(t)\right]\right.\\
      &&\times \left. \chi_1(2^{j} t - \varrho) \right\rvert \\
    \end{eqnarray*} 
\end{theorem}
 Again, we have from the from Hardy–Littlewood–Sobolev inequality \cite{In1} and by Theorem \eqref{TH1}, we have
\begin{eqnarray*}
    \lvert ( \mathfrak{W}_{\varrho}^{\varkappa_1,\chi_2} h)(t) - h(t) \rvert &\leq & 2^{\varrho -1}\digamma_{\upsilon}^{-1} G(p,\varrho)  |t_1-t_2|^{-\varrho}\left[\lVert \upsilon(t) t^{\varrho-1} \rVert_1\lVert \theta(t) \rVert_1 +\lVert \theta(t) t^{\varrho-1} \rVert_1\lVert \lVert \upsilon(w) \rVert_1\right]\Upsilon(w)\rVert_1\\
   &\times&\sum_{\varrho =-\infty}^{+\infty} \left[ \chi_1(2^{j} t - k)\int_{\mathbb{R}} |t-u|^{-\varrho} \left. \chi_2(2^{j} u - \varrho) \, du\int_{\infty} f(t_3) \, dt_3 \right] - 2^{-\frac{J}{2}} h(t) \right|\\
   &\leq & 2^{\varrho -1}\digamma_{\upsilon}^{-1} G(p,\varrho) \left[\lVert \upsilon(t) t^{\varrho-1} \rVert_1\lVert \theta(t) \rVert_1 +\lVert \theta(t) t^{\varrho-1} \rVert_1\lVert \lVert \upsilon(w) \rVert_1\right]\Upsilon(w)\rVert_1\\
   &\times& \left| \sum_{\varrho =-\infty}^{+\infty} \left[ \int_{\mathbb{R}} \int_{\mathbb{R}} \frac{f(t_3)\chi_2(2^{j} u - \varrho)}{|t-u|^{\varrho}} \, du \, dt_3 \right] \right| \lvert \varkappa_1(2^{j} t - \varrho) \rvert \, (2^{\frac{\varrho}{2}}+1)\\
  &\leq & 2^{\varrho -1}\digamma_{\upsilon}^{-1} G^{'}(p,\varrho) \left[\lVert \upsilon(t) t^{\varrho-1} \rVert_1\lVert \theta(t) \rVert_1 +\lVert \theta(t) t^{\varrho-1} \rVert_1\lVert \lVert \upsilon(w) \rVert_1\right]\Upsilon(w)\rVert_1\\
  &\times& \lvert f(t_3) \rVert_r \lVert \chi_2 \rVert_{r'},
\end{eqnarray*} where $G^{'}(p,\varrho) $ is a real constant.

\end{theorem}
\begin{theorem}
    Let $\varkappa_1$\,\,$\chi_2$ be the wavelets having $m-$vanishing moments, compact support and $\varkappa_1,\,\chi_2\in \mathbb{C}^r$ for $r\leq m $. Further, assume that $h(t)\in \mathbb{C}^m$ such that for $K>0, \, \lvert h^m(t)\rvert \leq K$. The operators defined as \eqref{E4}, we have an approximation error as follows:
    \begin{equation}
        \label{Inq2}
        \lvert \mathfrak{W}_{\varrho}^{\varkappa_1, \chi_2} (h)(t)-h(t)\rvert \leq C2^{-jm}, 
    \end{equation}
    where $C=\frac{KC_1C_2}{m!}$ with the following meaning as;
    \begin{itemize}
         \item $C_1$ denotes $\int_{\mathbb{R}} |u|^m |\chi_2(u)|du$ is $m_{th}$ moment of $\chi_2$
          \item $C_2=\sum_{\varrho}\varkappa_1(2^j-\varrho)$ denotes the absolute sum which runs over the compact support.
    \end{itemize}
    \proof Given $h(t)\in \mathbb{C}^r$. Therefore, from Taylor's series of expansion, we have
    \begin{eqnarray*}
        \int_{\mathbb{R}} h(u)\chi_2(2^ju-\varrho)du&=& \sum_{k=0}^{m-1}\frac{h^{(k)}(t)}{k!}\int_{\mathbb{R}}(u-t)^k \chi_2(2^ju-k)du+\int_{\mathbb{R}} E_m(u) \chi_2(2^ju-k)du
    \end{eqnarray*}
    Since $\chi_2$ has $m$ vanishing moments i.e.,
$$\int_{\mathbb{R}}(u-t)^k \chi_2(2^ju-k)du=0,\,\, k=0,1,2,\cdots, m-1.$$
Thus the error is estimated as 
    \begin{eqnarray*}
        \lvert \mathfrak{W}_{\varrho}^{\chi_1, \chi_2} (h)(t)-h(t)\rvert\leq \left\lvert \sum_{\varrho=-\infty}^{\infty}\varkappa_1(2^ju-\varrho)\int_{\mathbb{R}} E_m(u) \chi_2(2^ju-k)du\right\rvert.
    \end{eqnarray*}
Obviously, $ E_m(u)\leq \left\lvert\frac{K (u-t)^m}{m!}\right\rvert$, we have
\begin{eqnarray*}
            \lvert \mathfrak{W}_{\varrho}^{\varkappa_1, \chi_2} (h)(t)-h(t)\rvert&\leq& \frac{K}{m!}\left\lvert \sum_{\varrho=-\infty}^{\infty}\chi_1(2^ju-\varrho)\int_{\mathbb{R}} (u-t)^m\chi_2(2^ju-k)du\right\rvert.
\end{eqnarray*}
In addition, assume that $2^ju-\varrho=v$. Then the above inequality reduces to
\begin{eqnarray*}
   \int_{\mathbb{R}} |u-t|^m|\chi_2(2^ju-k)|du\leq  \int_{\mathbb{R}} \left\lvert\left(\frac{v+\varrho}{2^j}-t\right)\right\rvert^m|\chi_2(v)|du
\end{eqnarray*}
The limit of integration is finite since the $\chi_2(t)$ has a compact support. Suppose $$C_1=\int_{\mathbb{R}}|v|^m \chi_2(v) dv,\,\, C_2=\sum_{\varrho=-\infty}^{\infty}\chi_1(2^ju-\varrho).$$ 
Now, \begin{eqnarray*}
   \int_{\mathbb{R}} |u-t|^m|\chi_2(2^ju-k)|du\leq  \frac{C_1}{2^{jm}}.
\end{eqnarray*}
So, \begin{eqnarray*}
            \lvert \mathfrak{W}_{\varrho}^{\chi_1, \chi_2} (h)(t)-h(t)\rvert&\leq& \frac{KC_1C_2}{m!}2^{-jm}.
\end{eqnarray*}
Thus we have
\begin{eqnarray*}
            \lvert \mathfrak{W}_{\varrho}^{\chi_1, \chi_2} (h)(t)-h(t)\rvert&\leq&C2^{-jm},
\end{eqnarray*}
where $$C=\frac{KC_1C_2}{m!}.$$
\end{theorem}
\textbf{Note 5} From the expression defined in \eqref{Inq2}, the dependency of $j$ shows that the approximation error decreases exponentially as $j$ increases. Whereas the dependency on $m$ reflects the error bound improves with increasing $m$.\par
\textbf{Note 6} From the practical point of view, if $h(t)$ is a smooth function, then the operators $\mathfrak{W}_{\varrho}^{\chi_1, \chi_2} (h)(t)$ provide a small error approximation, but near sharp edges (discontinuities) of $h(t)$, the error is unpredictable.
\subsection{Error Estimation of the Sampling Theorem}
For $ \mathfrak{W}_{\varrho}^{\varkappa_1,\chi_2} h  \in V_J(\chi_1)$, the operators defined in \eqref{E4} do not always work. However, we can still construct the following operators \( (\mathfrak{W}_{\varrho}^{\varkappa_1,\chi_2} h )(t) \) from the samples \( (\mathfrak{W}_{\varrho}^{\varkappa_1,\chi_2} h )\left(\frac{\varrho}{2^J}\right) \) as defined in \eqref{E4}. 
Furthermore, if \( \mathfrak{W}_{\varrho}^{\varkappa_1,\chi_2} \notin V_{J+1}(\chi_1) \) and \( \mathfrak{W}_{\varrho}^{\varkappa_1,\chi_2}   \in V_J(\chi_1) \),

\begin{eqnarray}\label{samp}
    \Tilde{(\mathfrak{W}}_{\varrho}^{\chi_1,\chi_2} h )(t) =  \sum_{\varrho=-\infty}^{+\infty} \varkappa_1(2^{j} t - \varrho) \int_{\mathbb{R}} \Tilde{h}(u) \chi_2(2^{j} u - \varrho) \, du
\end{eqnarray}
Walter \cite{W1} considered the error estimate
$\lVert \Tilde{\mathfrak{W}}_{\varrho}^{\varkappa_1,\chi_2}( h) - \mathfrak{W}_{\varrho}^{\varkappa_1,\chi_2} (h) \rVert$  \text{for more general} \quad $\mathfrak{W}_{\varrho}^{\varkappa_1,\chi_2}\\\in L^2(\mathbb{R}) \quad \text{than} \quad \mathfrak{W}_{\varrho}^{\varkappa_1,\chi_2} \in V_J(\chi).$
Here, the scaling function $\chi(t)$ may not necessarily be a cardinal orthogonal spline function (COSF). To handle this, we first decompose the space $L^2(\mathbb{R})$ as follows:

For $\delta>0$,
\begin{eqnarray}
    U_{J}(\delta)\underline{\approx} \{ f\in L^2(\mathbb{R});(1-\chi_{[-\pi,\pi]}(\omega)\hat{f}(2^{j}\omega)\}\leq \delta.
\end{eqnarray}
In that order, we have the following proposition in terms of $B_1(f,\chi)$\\
\textbf{Proposition 2} For $\delta>0$, and $\mathfrak{W}_{\varrho}^{\varkappa_1,\chi_2} \in U_J({\delta})$, we have
\begin{eqnarray*}
    B_1(\mathfrak{W}_{\varrho}^{\varkappa_1,\chi_2} (h),\chi)-\frac{2}{\sqrt{2\pi}}\delta\leq \lVert \Tilde{\mathfrak{W}}_{\varrho}^{\varkappa_1,\chi_2}( h) - \mathfrak{W}_{\varrho}^{\varkappa_1,\chi_2} (h) ) \rVert\leq B_1(\mathfrak{W}_{\varrho}^{\varkappa_1,\chi_2} (h),\chi)+\frac{2}{\sqrt{2\pi}}\delta,
\end{eqnarray*} where $[{\mathfrak{W}_{\varrho}^{\varkappa_1,\chi_2} (h)}]_J(t)$ is orthogonal projection of $\mathfrak{W}_{\varrho}^{\varkappa_1,\chi_2} (h)$ on $V_j(\chi)$ and \begin{equation}\label{In1}
     B_1(f,\chi)=\sqrt{\frac{2^{J-1}}{\pi}}\left(\int_{-\pi}^{\pi} \hat{\mathfrak{W}}_{\varrho}^{\chi_1,\chi_2} \lvert(-2^{j}\omega))\rvert^2 |1-\Tilde{\chi}(\omega^2|^2 d\omega\right)^{\frac{1}{2}}.
\end{equation}
\label{L1}\textbf{Lemma 1} For $\varrho \in \mathbb{Z}$,
$$c_{j\varrho}-2^{\frac{-J}{2}}f\left( \frac{\varrho}{2^J}\right)=\frac{2^{\frac{J}{2}}}{2\pi}\int_{\mathbb{R}} \hat{\mathfrak{W}}_{\varrho}^{\chi_1,\chi_2}(-2^{j}\omega)(\chi(\omega)-1)e^{-i\varrho \omega}d\omega,$$
where $c_{j\varrho}$ is defined \eqref{4}.
\proof
We have,
\begin{eqnarray*}
    c_{jk}&=&2^{\frac{J}{2}}\int_{\mathbb{R}} \hat{\mathfrak{W}}_{\varrho}^{\chi_1,\chi_2}(t)\chi(2^{j}t-\varrho) d t \\
   &=& \frac{2^{\frac{J}{2}}}{2\pi}\int_{\mathbb{R}} \hat{\mathfrak{W}}_{\varrho}^{\chi_1,\chi_2}(-2^{j}\omega)(\chi(\omega)-1)e^{-i\varrho \omega}d\omega
\end{eqnarray*}
Again, 
\begin{eqnarray*}
    2^{\frac{-J}{2}}f\left( \frac{\varrho}{2^J}\right) &=& \frac{2^{\frac{J}{2}}}{2\pi}\int_{\mathbb{R}} \hat{\mathfrak{W}}_{\varrho}^{\chi_1,\chi_2}(-2^{j}\omega)\chi(\omega)e^{-i\varrho \omega}d\omega
\end{eqnarray*}
\begin{theorem} \label{err}
     The scaling function $\chi_1(t)$ is the operator defined as in \eqref{E4}. Then
     \begin{eqnarray*}
         \lVert \Tilde{\mathfrak{W}}_{\varrho}^{\varkappa_1,\chi_2} (h) - \mathfrak{W}_{\varrho}^{\varkappa_1,\chi_2} (h ) \rVert^2= \lVert \mathfrak{W}_{\varrho}^{\varkappa_1,\chi_2} (h) - [\mathfrak{W}_{\varrho}^{\varkappa_1,\chi_2} (h )]_J \rVert^2 + B_2( \Tilde{\mathfrak{W}}_{\varrho}^{\varkappa_1,\chi_2} (h) - \mathfrak{W}_{\varrho}^{\varkappa_1,\chi_2} (h ),\phi))^2,
     \end{eqnarray*}
     where $\lVert \mathfrak{W}_{\varrho}^{\varkappa_1,\chi_2} (h) - [\mathfrak{W}_{\varrho}^{\varkappa_1,\chi_2} (h )]_J \rVert$ the above proposition 
     
     \begin{eqnarray*}
    B_2(\mathfrak{W}_{\varrho}^{\varkappa_1,\chi_2} (h ),\phi) = \frac{2^{\frac{2}{J}}}{2\pi} ( \sum_{k} \left\lvert \int_{\mathbb{R}} \left( 1 - \Tilde{\chi}_1(\omega)) \Tilde{f}(-2^J \omega)  e^{-ik\omega} \, d\omega \right\rvert^2 \right)^{\frac{1}{2}}
\end{eqnarray*}

\proof From the \eqref{samp}
\begin{eqnarray}\label{samp}
    \Tilde{(\mathfrak{W}}_{\varrho}^{\chi_1,\chi_2} h )(t) = \sum_{\varrho=-\infty}^{+\infty} 2^{\frac{-J}{2}} [\chi_1]_{J\varrho}(t) \int_{\mathbb{R}} \Tilde{h}(u) \chi_2(2^{j} u - \varrho) \, du,
\end{eqnarray}
where $[\chi_1]_{J\varrho}(t)= 2^{\frac{-J}{2}} \chi_1(2^{j}t-\varrho)$. From \eqref{4} and Lemma 1, we have
\begin{eqnarray*}
    \lVert \Tilde{(\mathfrak{W}_{\varrho}^{\varkappa_1,\chi_2}}-\mathfrak{W}_{\varrho}^{\varkappa_1,\chi_2}\rVert&=&\left(\lVert \sum_{\varrho=-\infty}^{+\infty} 2^{\frac{-J}{2}}  \int_{\mathbb{R}} \Tilde{h}(u) \chi_2(2^{j} u - \varrho) \, du-[\chi_1]_{J\varrho}(t)\rVert ^2\right)^{\frac{1}{2}}\\
    &=& \frac{2^{\frac{2}{J}}}{2\pi}  \sum_{k} \left\lvert \int_{\mathbb{R}} \left( \left(1 - \Tilde{\chi}_1(\omega)\right) \Tilde{f}(-2^J \omega)  e^{-ik\omega} \, d\omega \right\rvert^2 \right)^{\frac{1}{2}}.
\end{eqnarray*}

\end{theorem}

\textbf{Proposition 3} \label{errcc}
    If,\( \mathfrak{W}_{\varrho}^{\varkappa_1,\chi_2} \in U_{J}(\delta) \) with \( \delta = 0 \), or equivalently, if \( \mathfrak{W}_{\varrho}^{\varkappa_1,\chi_2}( h) \) is \( 2^{j}\pi \) band-limited \,\&\, \( \chi_1 \) is a scaling function, then
    \begin{eqnarray*}
        \lVert \Tilde{\mathfrak{W}}_{\varrho}^{\varkappa_1,\chi_2}( h) - \mathfrak{W}_{\varrho}^{\varkappa_1,\chi_2} (h) ) \rVert=B_1(\mathfrak{W}_{\varrho}^{\varkappa_1,\chi_2} (h),\chi_1)^2+B_3(\mathfrak{W}_{\varrho}^{\varkappa_1,\chi_2} (h),\chi_1)^2\\
    \end{eqnarray*}
 where 
    \begin{eqnarray*}
    B_3(\mathfrak{W}_{\varrho}^{\varkappa_1,\chi_2}( h) ,\chi)=\sqrt{\frac{2^{J-1}}{\pi}}\left(\int_{-\pi}^{\pi} \hat{\mathfrak{W}}_{\varrho}^{\chi_1,\chi_2} \lvert(-2^{j}\omega))\rvert^2 |1-\Tilde{\chi}(\omega^2|^2 d\omega\right)^{\frac{1}{2}}.
\end{eqnarray*}
\proof Given that if \( \mathfrak{W}_{\varrho}^{\varkappa_1,\chi_2}\in U_{J}(\delta) \) with \( \delta = 0 \),
\begin{eqnarray*}
     \sum_{k} \left\lvert \int_{\mathbb{R}}  \left(1 - \Tilde{\chi}_1(\omega)\right) \Tilde{ \mathfrak{W}}_{\varrho}^{\chi_1,\chi_2}( h) (-2^J \omega)  e^{-ik\omega} \, d\omega \right\rvert^2 
     &=&\sum_{k} \left\lvert \int_{-\pi}^{{\pi} }  \left(1 - \Tilde{\chi}_1(\omega)\right) \Tilde{ \mathfrak{W}}_{\varrho}^{\chi_1,\chi_2}( h) (-2^J \omega)  e^{-ik\omega} \, d\omega \right\rvert^2 \\
     &=& 2\pi \int_{-\pi}^{{\pi} }\left\lvert  \left(1 - \Tilde{\chi}_1(\omega)\right) \Tilde{ \mathfrak{W}}_{\varrho}^{\chi_1,\chi_2}( h) (-2^J \omega)  e^{-ik\omega}  \right\rvert^2\, d\omega \\
\end{eqnarray*}
   Hence, using the Theorem \eqref{err}, we have the required result.
\section{Illustration:}\label{EXA}
\begin{example}
    Let us define the second-order B-spline $h\in \mathbb{C}^{\mathbb{N}}(\mathbb{R})$ as follows: \[
   h (t) = 
   \begin{cases} 
   \frac{1}{2} t^2, & \text{if } 0 \leq t < 1, \\
   \frac{1}{2} (2 - t)^2, & \text{if } 1 \leq t < 2, \\
   0, & \text{otherwise.}
   \end{cases}
   \] Clearly, h(t) is compactly supported up to the second order. Then the first term, namely $\sum_{i=1}^{\varrho} \frac{\lvert h^i(x)\rvert}{i!}\frac{a^i}{2^{i\varrho}}$ in the inequality \eqref{T1} elaborates the exponential decay factor $\frac{a^i}{2^{i\varrho}}$ that ensures that the series converges rapidly for a large value of $\varrho$. Again noting the fact that $h(t), h'(t) \, \& h^{''}(t)$ are bounded functions and the second term $ \frac{a^\mathbb{N}}{2^{\mathbb{N}\varrho}}\omega_1\left( h^{\mathbb{N}}, \frac{a}{2^{j}}\right)$ in the inequality \eqref{T1} indicates that the modulus of smoothness depends on the smoothness of h(t). Particularly, for the B-spline of order 2, this term is small because of the boundedness of higher derivatives and the compact support.\par
  Clearly, from the basic reasoning one can easily come to the conclusion of verifying the fundamental theorem of approximation in terms of the operators point of view using this B-spline function. In other words, both terms, namely in RHS $\sum_{i=1}^{\varrho} \frac{\lvert h^i(x)\rvert}{i!}\frac{a^i}{2^{i\varrho}}+ \frac{a^\mathbb{N}}{2^{\mathbb{N}\varrho}}\omega_1\left( h^{\mathbb{N}}, \frac{a}{2^{j}}\right)$, will tend to zero for a sufficiently large value of $\varrho$. That is $  \mathfrak{W}_{\varrho}^{\varkappa_1,\chi_2}h \rightarrow h \, \forall t\in \mathbb{R}$, for the long run of $\varrho$.
\end{example}
\begin{example}
Suppose that cubic B-spline 
  
   \[
f(t) =
\begin{cases} 
\frac{1}{6} (2 - |t|)^3, & \text{if } 0 \leq t \leq 2, \\
0, & \text{otherwise.}
\end{cases}
\]

The Fourier Transform of the B-spline function is given by the following expression:

\[
\tilde{f}(\omega) = \left(\frac{\sin\frac{\omega}{2}}{\frac{\omega}{2}}\right)^{4}
\]
Particularly, for the $J=2$ and $\delta=0.01$, let us find lower and upper bound of the following inequality 
   \begin{eqnarray*}
        B_1(\mathfrak{W}_{\varrho}^{\varkappa_1,\chi_2} (h),\chi)-\frac{2}{\sqrt{2\pi}}\delta\leq \lVert \Tilde{\mathfrak{W}}_{\varrho}^{\varkappa_1,\chi_2}( h) - \mathfrak{W}_{\varrho}^{\varkappa_1,\chi_2} (h) ) \rVert\leq B_1(\mathfrak{W}_{\varrho}^{\varkappa_1,\chi_2} (h),\chi)+\frac{2}{\sqrt{2\pi}}\delta, 
   \end{eqnarray*}
   where $B_1(\mathfrak{W}_{\varrho}^{\varkappa_1,\chi_2} (h),\chi)$ is defined as in \eqref{In1}. The above expression takes a new form after doing the significant calculation over Python as;
   \begin{eqnarray*}
       0.00866839702279679 \leq \lVert \Tilde{\mathfrak{W}}_{\varrho}^{\varkappa_1,\chi_2}( h) - \mathfrak{W}_{\varrho}^{\varkappa_1,\chi_2} (h) ) \rVert\leq 0.024626088238854097. 
   \end{eqnarray*}
\end{example}
\begin{example}
    Assume that $h(t)=\sin(2\pi t),\,\,t\in[0,1]$. Now, we observe the nature of the operators $\Tilde{\mathfrak{W}}_{\varrho}^{\varkappa_1,\chi_2}( h)$ approximate to the function $h(t)$ near the sharp edges as seen in the following table and graph by denoting error as $E(t)$, which is given as $E(t)=\lvert(\Tilde{\mathfrak{W}}_{\varrho}^{\varkappa_1,\chi_2} h)(t)-h(t)\rvert$

\setlength{\arrayrulewidth}{0.5mm}
\setlength{\tabcolsep}{18pt}
\renewcommand{\arraystretch}{2}

\begin{table}[h]
    \centering
    \renewcommand{\arraystretch}{1.8} 
    \begin{tabular}{|c|c|c|c|c|c|}
        \hline
        \multicolumn{6}{|c|}{\textbf{Value of $E(t)=\lvert(\Tilde{\mathfrak{W}}_{\varrho}^{\varkappa_1,\chi_2} h)(t)-h(t)\rvert$ for following $t$s and $\varrho$s }} \\
        \hline
        $t$ & $\varrho = 0$ & $\varrho = 1$ & $\varrho = 2$ & $\varrho = 3$ & $\varrho = 4$ \\
        \hline
        0.1000 & 0.6537 & 0.5878 & 0.5878 & 0.5878 & 0.5878 \\
        \hline
        0.3000 & 0.9511 & 0.8851 & 0.9511 & 0.9511 & 0.9511 \\
        \hline
        0.5000 & 0.0000 & 0.0000 & 0.0659 & 0.0000 & 0.0000 \\
        \hline
        0.7000 & 0.9511 & 0.9511 & 0.8851 & 0.9511 & 0.9511 \\
        \hline
        0.9000 & 0.5878 & 0.5878 & 0.5878 & 0.5878 & 0.5878 \\
        \hline
    \end{tabular}
    \caption{ E(t) with respect to the function $f(t)= \sin(2\pi t)$}
    \label{tab:error_table}
\end{table}

\begin{figure}
    \centering
    \includegraphics[width=1\textwidth]{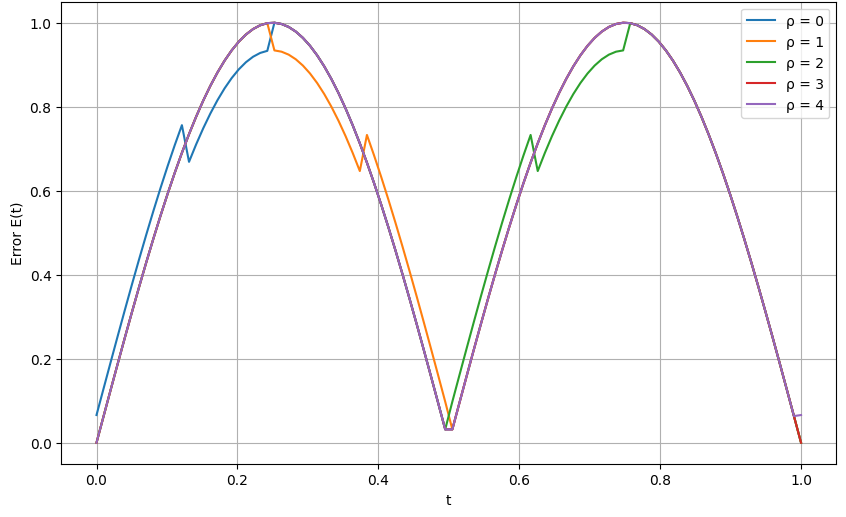}

    \caption{Graphical representation of E(t) with respect to the function $h(t)= \sin(2\pi t)$.}
    \label{fig:fejer_results}
\end{figure}

\end{example}
\newpage

\newpage
We summarize from Table $1$ and Figure $1$, that the error $E(t)$ for different $\varrho$s has an unorthodox behavior  at $t=0.5$. More preciously,  the operators $\Tilde{\mathfrak{W}}_{\varrho}^{\varkappa_1,\chi_2} f)(t)$ approximate to $f(t)$ better in the smooth region of $f(t)$ for increasing values of $\varrho$ with exception near the point of discontinuity.

\section{Conclusion} In this article, we have made an effort to construct wavelet-based filtering operators and establish the related approximation theorems with some remarks and important notes. We deeply explored the COSF for the wavelet-based sampling theorem in MRA and observed that standard wavelet sampling may not work if a signal does not come from the multiresolution spaces. Later, a theorem based on an error estimate along with lower and upper bound, is derived on the basis of the decomposition of the subspaces. In that order, we came to the fact that there is an increment in error estimate near the sharp edges of a function, whereas a better error estimate has been found on the smooth portion of a continuously and discontinuously differentiable function. 

\vspace{0.55cm}
\begin{center}
\textbf{Acknowledgments}
\end{center}
The first author expresses his heartfelt gratitude to the students and lab mates from the Department of Applied Sciences and Humanities, Institute of Engineering and Technology, Lucknow, Uttar Pradesh, India, for their continuous discussions and invaluable insights that kept him motivated throughout this research.

Furthermore, he acknowledges the financial support provided by the ``Homi Bhabha Teaching cum Research Fellowship" from Dr. A. P. J. Abdul Kalam Technical University, Lucknow, India.
\vspace{0.35cm}
\begin{center}
\textbf{Compliance with Ethical Standards}
\end{center}
\textbf{Funding:} No funding was received to report this study.\\
\textbf{Conflict of interest:} All the authors declare that no conflict of interest exists.\\
\textbf{Ethical approval:} This article does not contain any studies with human participants or animals performed by any of the authors.\\
\textbf{Data availability:} We assert that no data sets were generated or analyzed during the preparation of the manuscript.\\
\textbf{Code availability:} Not applicable.\\
\textbf{Authors' contributions:} All the authors have equally contributed to the conceptualization, framing, and writing of the manuscript.

\end{document}